\documentclass[reqno]{amsproc}
\usepackage{amssymb,amsmath,enumerate,mathrsfs}

\numberwithin{equation}{section}

\newtheorem{theorem}[equation]{Theorem}
\newtheorem{proposition}[equation]{Proposition}

\newtheorem{corollary}[equation]{Corollary}

\theoremstyle{definition}
\newtheorem{definition}[equation]{Definition}

\def\C{\mathbb C}
\def\HT{(\mathcal {HT})}
\def\L{\mathscr L}
\def\N{\mathbb N}
\def\R{\mathbb R}
\def\S{\mathscr S}
\def\x{\mathbf x}
\def\V{\mathcal V}
\def\Mbar{\overline{M}}
\def\RR{\mathcal R}

\DeclareMathOperator{\sym}{ \sigma\!\!\!\sigma}
\DeclareMathOperator{\End}{End}

\def\cu{\textup{cu}}
\def\cusym{\,{}^{\cu}\!\sym}
\def\cuT{\,{}^{\cu} T}
\def\cupi{\,{}^{\cu}\hspace{-1.5pt}\pi}
\def\cuH{{}^{\cu}H}
\def\cuPsi{{}^{\cu}\Psi}
\def\cl{\textup{cl}}
\def\eps{\varepsilon}

\def\st{;\;}

\DeclareMathOperator{\cuDiff}{{}^{\cu}Diff}

\DeclareMathOperator{\spec}{spec}
\DeclareMathOperator{\op}{op}

\begin{document}

\nocite*

\title[Maximal regularity for parabolic PDEs on mfds. with cyl. ends]{Maximal $L^p$--$L^q$ regularity
for parabolic partial differential equations on manifolds with cylindrical ends}
\author{Thomas Krainer}
\address{Mathematics and Statistics \\ Penn State Altoona \\
3000 Ivyside Park \\ Altoona, PA 16601 \\ U.S.A.}
\email{krainer@psu.edu}

\begin{abstract}
We give a short, simple proof of maximal $L^p$--$L^q$ regularity for linear parabolic 
evolution equations on manifolds with cylindrical ends by making use of
pseudodifferential parametrices and the concept of $\RR$-boundedness for
the resolvent.
\end{abstract}


\maketitle

\section{Introduction and main results}

\noindent
Let $(M,g)$ be a Riemannian manifold with cylindrical ends, i.e., there exists
a relatively compact, open subset $K \subset M$ such that $M \setminus K$ is isometric to
$[0,\infty)_r\times Y$, where $(Y,g_Y)$ is a closed compact Riemannian manifold (not
necessarily connected), and the cylinder $[0,\infty)_r\times Y$ is equipped with the
product metric $dr^2 + g_Y$. By employing the change of variables $x = 1/r$ on
$[0,\infty)_r\times Y$ for large values of $r$ and attaching a copy of $Y$
at $x = 0$, we obtain a compactification of $M$ to a smooth compact manifold
$\Mbar$ with boundary $\partial\Mbar = Y$. The interior of $\Mbar$ is diffeomorphic
to the original manifold $M$, and in a collar neighborhood $[0,\eps)_x\times Y$ of
the boundary the Riemannian metric now takes the form
${}^\cu g = \frac{dx^2}{x^4} + g_Y$, which is the form of a cusp metric (see \cite{Melrose}).
Differential operators on a manifold with cylindrical ends with a
`reasonably nice' coefficient behavior at infinity correspond in this
way to cusp differential operators on a compact manifold with boundary.
Similarly, function spaces on the original manifold $(M,g)$ with cylindrical
ends correspond to function spaces on $(\Mbar,{}^\cu g)$. Let $E \to \Mbar$ be a smooth
vector bundle. We will prove the following result:

\begin{theorem}\label{Maintheorem}
Let $A_t \in C([0,T],\cuDiff^m(\Mbar,E))$, $m > 0$, $ 0 < T < \infty$, and assume that
$A_t$ is  cusp-elliptic with parameter in $\Lambda = \{\lambda \in \C \st
\Re(\lambda) \geq 0\}$ for every $0 \leq t \leq T$. Then, for every $1 < p,q < \infty$,
the parabolic evolution equation
\begin{equation}\label{ParabEvolEqu}
\left.
\begin{aligned}
\frac{\partial u}{\partial t} - A_tu &= f \\
u|_{t=0} &= u_0
\end{aligned}
\right\}
\end{equation}
has a unique solution
\begin{align*}
u &\in W^{1,p}([0,T],{}^\cu L^q(\Mbar,E))\cap L^p([0,T],\cuH^{m,q}(\Mbar,E)) \\
\intertext{for every}
f &\in L^p([0,T],{}^\cu L^q(\Mbar,E)) \quad \textup{and} \\
u_0 &\in {}^\cu B_{q,p}^{m(1-1/p)}(\Mbar,E) = ({}^\cu L^q(\Mbar,E),\cuH^{m,q}(\Mbar,E))_{1-1/p,p}.
\end{align*}
The mapping $u \mapsto (f,u_0)$ given by \eqref{ParabEvolEqu} is a topological isomorphism
of these function spaces. In particular, the solution $u$ satisfies optimal $L^p$--$L^q$
a priori estimates with respect to the data $f$ and $u_0$.
\end{theorem}

By general results from \cite{PruessSchnaubelt} (see also \cite{Arendt} for recent
improvements), we need to prove Theorem~\ref{Maintheorem} only in the autonomous
case $A_t \equiv A$ and for $u_0 = 0$, i.e., we need to prove that if
$A \in \cuDiff^m(\Mbar,E)$, $m > 0$, is cusp-elliptic with parameter in
$\Lambda = \{\lambda \in \C \st \Re(\lambda) \geq 0\}$, then $A$ has maximal
$L^p$--$L^q$ regularity on $[0,T]$. By \cite{Weis} this is the case if the resolvent
$(A-\lambda)^{-1} : {}^{\cu}L^q(\Mbar,E) \to \cuH^{m,q}(\Mbar,E)$ exists for
$\lambda \in \Lambda$ with $|\lambda| \geq R$ sufficiently large, and if in addition
the family
\begin{equation}\label{RbddRes}
\{\lambda (A-\lambda)^{-1} \st \lambda \in \Lambda,\; |\lambda| \geq R\} \subset
\L\bigl({}^\cu L^q(\Mbar,E)\bigr)
\end{equation}
is $\RR$-bounded.
Thus Theorem~\ref{Maintheorem} is a consequence of the following result,
which is of interest in its own right.

\begin{theorem}\label{RbddnessResolvents}
Let $\Lambda \subset \C$ be a closed sector, and let $A \in \cuDiff^m(\Mbar,E)$, $m > 0$,
be cusp-elliptic with parameter in $\Lambda$. Let $1 < q < \infty$. For
$\lambda \in \Lambda$ with $|\lambda| \geq R$ sufficiently large, the operator
$$
A - \lambda : \cuH^{m,q}(\Mbar,E) \to {}^{\cu}L^q(\Mbar,E)
$$
is invertible, and the set \eqref{RbddRes} is $\RR$-bounded.
\end{theorem}

To prove Theorem~\ref{RbddnessResolvents} we will employ a parameter-dependent
parametrix of $A - \lambda$ in the calculus of cusp pseudodifferential operators
on $\Mbar$ to approximate the resolvent. The parametrix is then further analyzed making
use of the results from \cite{DenkKrainer} on $\RR$-boundedness of families of
pseudodifferential operators.

\medskip

\noindent
The structure of this paper is as follows:

In Section~\ref{RbddnessReview} we review the definition and necessary results
about $\RR$-boundedness of operator families that we need. For a comprehensive account on 
general aspects of $\RR$-boundedness and its applications to parabolic
equations we refer to the monograph \cite{DenkHieberPruess} or
the survey paper \cite{KunstmannWeis}, for $\RR$-boundedness of families of
pseudodifferential operators see \cite{DenkKrainer} (see also \cite{HytoenenPortal} for 
related work).

Section~\ref{Cuspoperators} is devoted to cusp differential and
pseudodifferential operators on manifolds with boundary
(see \cite{LauterMoroianu, MelroseNistor}). Finally, Section~\ref{ProofTheorem}
contains the proof of Theorem~\ref{RbddnessResolvents}.

\section{$\RR$-boundedness and families of pseudodifferential operators}
\label{RbddnessReview}

\begin{definition}\label{Rboundednessdef}
Let $E$ and $F$ be Banach spaces. A subset ${\mathcal T} \subset \L(E,F)$ is called
$\RR$-bounded, if for some constant $C \geq 0$
the inequality
\begin{equation}\label{Rboundedinequality}
\Bigl(\sum\limits_{\eps_1,\ldots,\eps_N \in \{-1,1\}}\Bigl\|\sum\limits_{j=1}^N\eps_j T_je_j\Bigr\|\Bigr)
\leq C \Bigl(\sum\limits_{\eps_1,\ldots,\eps_N \in \{-1,1\}}\Bigl\|\sum\limits_{j=1}^N\eps_j e_j\Bigr\|\Bigr)
\end{equation}
holds for all choices of $T_1,\ldots,T_N \in {\mathcal T}$ and $e_1,\ldots,e_N \in E$, $N \in \N$.

The best constant
\begin{gather*}
C = \sup\Bigl\{\Bigl(\sum\limits_{\eps_1,\ldots,\eps_N \in \{-1,1\}}\Bigl\|\sum\limits_{j=1}^N\eps_j T_je_j\Bigr\|\Bigr) \st
N \in \N,\: T_1,\ldots,T_N \in {\mathcal T}, \\
\Bigl(\sum\limits_{\eps_1,\ldots,\eps_N \in \{-1,1\}}\Bigl\|\sum\limits_{j=1}^N\eps_j e_j\Bigr\|\Bigr) = 1\Bigr\}
\end{gather*}
in \eqref{Rboundedinequality} is called the $\RR$-bound of ${\mathcal T}$ and will be
denoted by $\RR({\mathcal T})$.
\end{definition}

The general properties of $\RR$-bounded sets yield to the following result
about functions with $\RR$-bounded range (see \cite{DenkKrainer},
Propositions 2.11 and 2.13):

\begin{proposition}\label{linftyproj}
Let $\Gamma$ be a nonempty set. Define $\ell_{\RR}^{\infty}(\Gamma,\L(E,F))$ as the space
of all functions $f : \Gamma \to \L(E,F)$ with $\RR$-bounded range and norm
\begin{equation}\label{Rboundnorm}
\|f\|_{\ell^{\infty}_{\RR}}:= \RR\bigl(f(\Gamma)\bigr).
\end{equation}
Then $\bigl(\ell_{\RR}^{\infty}(\Gamma,\L(E,F)),\|\cdot\|_{\ell^{\infty}_{\RR}}\bigr)$
is a Banach space, and the embeddings
$$
\ell^{\infty}(\Gamma) \hat{\otimes}_{\pi} \L(E,F) \hookrightarrow
\ell^{\infty}_{\RR}(\Gamma,\L(E,F)) \hookrightarrow
\ell^{\infty}(\Gamma,\L(E,F))
$$
are well defined and continuous.

The norm in $\ell^{\infty}_{\RR}$ is submultiplicative, i.e.
$$
\|f\cdot g\|_{\ell^{\infty}_{\RR}} \leq \|f\|_{\ell^{\infty}_{\RR}} \cdot \|g\|_{\ell^{\infty}_{\RR}}
$$
whenever the composition $f\cdot g$ makes sense, and we have $\|{\mathfrak 1}\|_{\ell^{\infty}_{\RR}} = 1$
for the constant map ${\mathfrak 1} \equiv \textup{Id}_E$.
\end{proposition}

\begin{corollary}[Corollary 2.14 in \cite{DenkKrainer}]\label{BilderRbounded}
\begin{enumerate}[i)]
\item Let $M$ be a smooth manifold, and let $K \subset M$ a compact subset. Let $f \in C^{\infty}(M,\L(E,F))$.
Then the range $f(K)$ is an $\RR$-bounded subset of $\L(E,F)$.
\item Let $f \in \S(\R^n,\L(E,F))$. Then the range $f(\R^n) \subset \L(E,F)$ is $\RR$-bounded.
\end{enumerate}
\end{corollary}
\begin{proof}
The assertion follows from Proposition \ref{linftyproj} in view of
\begin{align*}
C^{\infty}(M,\L(E,F)) &\cong C^{\infty}(M) \hat{\otimes}_{\pi} \L(E,F), \\
\S(\R^n,\L(E,F)) &\cong \S(\R^n) \hat{\otimes}_{\pi} \L(E,F).
\end{align*}
\end{proof}

In what follows all Banach spaces are assumed to be of class $\HT$ and to satisfy
Pisier's property $(\alpha)$ (this is of relevance for the validity of
Theorem~\ref{IterationRboundedness} further below). We do not supply these
definitions here, but merely note that all Banach spaces that are isomorphic
to a scalar $L^q$-space for some $1 < q < \infty$ have both properties, and
so do vector valued $L^q$-spaces provided that the target space satisfies
both properties. This will be sufficient for our purposes.

\begin{definition}\label{symboldef}
Let $\ell \in \N$ be fixed. For every $\mu \in \R$ we define the
anisotropic $\RR$-bounded symbol class $S^{\mu;\ell}_{\RR}(\R^n\times\R^q;E,F)$
to consist of all operator functions $a \in C^{\infty}(\R^n\times\R^q,\L(E,F))$ such that
for all $\alpha \in \N_0^n$ and $\beta \in \N_0^q$
\begin{equation}\label{SymbolR}
(1 + |\zeta| + |\lambda|^{1/\ell})^{-\mu+|\alpha|+\ell|\beta|}
\bigl(\partial_{\zeta}^{\alpha}\partial_{\lambda}^{\beta}a(\zeta,\lambda)\bigr)
\in \ell^{\infty}_{\RR}(\R^n\times\R^q,\L(E,F)).
\end{equation}
We equip this symbol class with the Fr{\'e}chet topology whose seminorms
are the $\ell^{\infty}_{\RR}$-norms of the functions in \eqref{SymbolR}.
\end{definition}

These spaces have all the usual properties of symbol spaces
(see \cite{DenkKrainer}, Section 3). Moreover, we have
$$
S^{\mu;\ell}_{\cl}(\R^n\times\R^q;E,F) \hookrightarrow S^{\mu;\ell}_{\RR}(\R^n\times\R^q;E,F)
\hookrightarrow S^{\mu;\ell}(\R^n\times\R^q;E,F)
$$
with continuous embeddings, where $S^{\mu;\ell}(\R^n\times\R^q;E,F)$ denotes the
standard space of anisotropic operator valued symbols of order $\mu$ (replace
$\ell^{\infty}_{\RR}$ in \eqref{SymbolR} by $\ell^{\infty}$), and
$S^{\mu;\ell}_{\cl}(\R^n\times\R^q;E,F)$ is the subspace of anisotropic
classical symbols, i.e., those that admit an asymptotic expansion
$a \sim \sum_{j=0}^{\infty}a_j$ with $a_j(\varrho\zeta,\varrho^{\ell}\lambda) = 
\varrho^{\mu-j}a_j(\zeta,\lambda)$ for $|(\zeta,\lambda)| \geq 1$ and $\varrho \geq 1$.
Furthermore,
$$
\bigcap\limits_{\mu \in \R}S^{\mu;\ell}_{\RR}(\R^n\times\R^q;E,F) =
\S(\R^n\times\R^q,\L(E,F)).
$$
For what we have in mind, the parameter space $\R^q$ needs to be replaced
by a closed sector $\Lambda \subset \R^2$. As is customary, the symbol spaces
in this case consist by definition of the restrictions of symbols defined in the
full space, and we equip those spaces with the quotient topology.

\medskip

Now split $\R^n = \R^d \times \R^{n-d}$ in the (co-)variables $\zeta = (\eta,\xi)$, where
$1 \leq d \leq n$ (in the case $d = n$ the $\R^{n-d}$-factor just drops out), and
consider symbols
\begin{equation}\label{SymbolsinSpace}
a(y,\eta,\xi,\lambda) \in S^{0}_{\cl}(\R^d_y,
S^{\mu;\ell}_{\RR}(\R^d_{\eta}\times\R^{n-d}_{\xi}\times\Lambda;E,F)).
\end{equation}
With $a(y,\eta,\xi,\lambda)$ we associate the family of pseudodifferential operators
\begin{align*}
A(\xi,\lambda) &= \op_y(a)(\xi,\lambda) : \S(\R^d,E) \to \S(\R^d,F),
\intertext{where}
[\op_y(a)(\xi,\lambda)u](y) &= (2\pi)^{-d}\iint e^{i(y-y')\eta}
a(y,\eta,\xi,\lambda)u(y')\,dy'\,d\eta.
\end{align*}

The following is a consequence of Theorem 3.18 in \cite{DenkKrainer}.

\begin{theorem}\label{IterationRboundedness}
For $\nu \geq \mu$ the family of pseudodifferential operators $\op_y(a)(\xi,\lambda)$ 
extends by continuity to
$$
\op_y(a)(\xi,\lambda) : H^{s,q}(\R^d,E) \to H^{s-\nu,q}(\R^d,F)
$$
for every $s \in \R$ and $1 < q < \infty$, and the operator function
$$
\R^{n-d}\times\Lambda \ni (\xi,\lambda) \mapsto \op_y(a)(\xi,\lambda) \in
\L\bigl(H^{s,q}(\R^d,E),H^{s-\nu,q}(\R^d,F)\bigr)
$$
belongs to the $\RR$-bounded symbol space $S_{\RR}^{\mu';\ell}(\R^{n-d}\times\Lambda;
H^{s,q}(\R^d,E),H^{s-\nu,q}(\R^d,F))$ with $\mu' = \mu$ if $\nu \geq 0$, or
$\mu' = \mu-\nu$ if $\nu < 0$.

The mapping $\op_y : a(y,\eta,\xi,\lambda) \mapsto \op_y(a)(\xi,\lambda)$ is continuous
in
$$
S^{0}_{\cl}(\R^d,S^{\mu;\ell}_{\RR}(\R^d\times\R^{n-d}\times\Lambda;E,F))
\!\to\!
S_{\RR}^{\mu';\ell}(\R^{n-d}\times\Lambda;H^{s,q}(\R^d,E),H^{s-\nu,q}(\R^d,F)).
$$
\end{theorem}

Let $Y$ be a closed compact manifold, and let $E,F \to Y$ be smooth (finite dimensional)
vector bundles. Let $L^{\mu;\ell}(Y,\R^{n-d}\times\Lambda;E,F)$ be the class of families
of pseudodifferential operators
$$
A(\xi,\lambda) : C^{\infty}(Y,E) \to C^{\infty}(Y,F)
$$
that are locally modelled on symbols \eqref{SymbolsinSpace}, and global remainders
on $Y$ that are integral operators with $C^{\infty}$-kernels that depend rapidly decreasing
(together with all derivatives) on the parameters
$(\xi,\lambda) \in \R^{n-d}\times\Lambda$. We write
$L^{\mu;\ell}_{\cl}(Y,\R^{n-d}\times\Lambda;E,F)$ if the symbols in \eqref{SymbolsinSpace}
are in addition required to be classical.
The following is an immediate consequence of Corollary~\ref{BilderRbounded} (applied
to global remainders) and Theorem~\ref{IterationRboundedness} (applied to families
supported in a chart).

\begin{corollary}\label{PseudosClosedRbdd}
Let $A(\xi,\lambda) \in L^{\mu;\ell}(Y,\R^{n-d}\times\Lambda;E,F)$. For $\nu \geq \mu$
and every $s \in \R$, $1 < q < \infty$, the operator family
$$
A(\xi,\lambda) : H^{s,q}(Y,E) \to H^{s-\nu,q}(Y,F)
$$
is continuous, and the operator function
$$
\R^{n-d}\times\Lambda \ni (\xi,\lambda) \mapsto A(\xi,\lambda) \in
\L\bigl(H^{s,q}(Y,E),H^{s-\nu,q}(Y,F)\bigr)
$$
belongs to the $\RR$-bounded symbol space $S_{\RR}^{\mu';\ell}(\R^{n-d}\times\Lambda;
H^{s,q}(Y,E),H^{s-\nu,q}(Y,F))$ with $\mu' = \mu$ if $\nu \geq 0$, or
$\mu' = \mu-\nu$ if $\nu < 0$.
The embedding
$$
L^{\mu;\ell}(Y,\R^{n-d}\times\Lambda;E,F) \hookrightarrow
S_{\RR}^{\mu';\ell}(\R^{n-d}\times\Lambda;H^{s,q}(Y,E),H^{s-\nu,q}(Y,F))
$$
is continuous.
\end{corollary}

\section{Analysis of cusp operators on manifolds with boundary}\label{Cuspoperators}

Let $\Mbar$ be a smooth $n$-dimensional compact manifold with boundary.
Let $U \cong [0,\eps)\times Y$ be a collar neighborhood of the boundary
$Y = \partial\Mbar$, and fix a smooth defining
function $\x$ for $Y$ (i.e. $\x \geq 0$ on $\Mbar$, $\x = 0$ precisely on $Y$, and
$d\x \neq 0$ on $Y$) that coincides in $U$ with the projection to the first
coordinate. With these choices we let
$$
{}^{\cu}\V = \{V \in C^{\infty}(\Mbar,T\Mbar) \st V\x \in \x^2C^{\infty}(\Mbar)\}
$$
be the Lie algebra of cusp vector fields on $\Mbar$. Let $\cuDiff^{\ast}(\Mbar)$ be
the envoloping algebra of cusp differential operators generated by ${}^{\cu}\V$ and
$C^{\infty}(\Mbar)$. In coordinates near the boundary, an operator
$A \in \cuDiff^{m}(\Mbar)$ takes the form
\begin{equation}\label{CuspDiffOp}
A = \sum\limits_{k+|\alpha| \leq m} a_{k,\alpha}(x,y)(x^2D_x)^kD_y^{\alpha}
\end{equation}
with $C^{\infty}$-coefficients $a_{k,\alpha}$ that are smooth up to $x = 0$.

The vector fields ${}^{\cu}\V$ are a projective finitely generated module over
$C^{\infty}(\Mbar)$, hence by Swan's theorem there is a smooth vector bundle
${}^{\cu}T\Mbar$ on $\Mbar$, the cusp-tangent bundle, whose space of
$C^{\infty}$-sections is ${}^{\cu}\V$. Locally near the boundary, the vector fields
$x^2\partial_x$ and $\partial_{y_j}$, $j = 1,\ldots,n-1$, form a frame for this bundle.
There is a canonical homomorphism $\phi : {}^{\cu}T\Mbar \to T\Mbar$ that restricts to
an isomorphism over the interior of $\Mbar$.
Let ${}^{\cu}T^*\Mbar$ be the cusp-cotangent bundle, the dual bundle to
${}^{\cu}T\Mbar$, and let ${}^t\phi : T^*\Mbar \to {}^{\cu}T^*\Mbar$ be the dual map
to $\phi$. Let $\sym(A)$ be the principal symbol of $A \in \cuDiff^m(\Mbar)$,
defined on $T^*\Mbar \setminus 0$. Over the interior of $\Mbar$ we set
$$
\cusym(A) = \sym(A) \circ \bigl({}^t\phi\bigr)^{-1}.
$$
This function extends by continuity to a smooth function on all of
${}^{\cu}T^*\Mbar \setminus 0$, and it is homogeneous of degree $m$ in the fibres.
$\cusym(A)$ is called the cusp-principal symbol of $A$. In coordinates near the
boundary, the cusp-principal symbol of the operator $A$ in \eqref{CuspDiffOp} is given
by
$$
\cusym(A) = \sum\limits_{k+|\alpha| = m} a_{k,\alpha}(x,y)\xi^k\eta^{\alpha}.
$$
More generally, if $E$ is a smooth vector bundle on $\Mbar$, we let $\cuDiff^m(\Mbar;E)$
denote the space of cusp differential operators $A$ of order (at most) $m$. Initially, we
consider $A$ an operator
\begin{equation}\label{AInitially}
A : \dot{C}^{\infty}(\Mbar,E) \to \dot{C}^{\infty}(\Mbar,E),
\end{equation}
where $\dot{C}^{\infty}(\Mbar,E)$ denotes the space of smooth sections of $E$ that vanish
to infinite order on the boundary.
The cusp-principal symbol of $A$ is a section $\cusym(A) \in
C^{\infty}\bigl(\cuT^*\Mbar\setminus 0,\End(\cupi^*E)\bigr)$, where
$\cupi : \cuT^*\Mbar\setminus 0 \to \Mbar$ is the canonical projection.

\begin{definition}\label{ElliptwithParam}
Let $\Lambda \subset \C$ be a closed sector. We call $A \in \cuDiff^m(\Mbar;E)$
cusp-elliptic with parameter in $\Lambda$ if
$$
\spec(\cusym(A))\cap\Lambda = \emptyset
$$
everywhere on $\cuT^*\Mbar\setminus 0$.
\end{definition}

\noindent
The operator \eqref{AInitially} extends by continuity to an operator
$$
A : \cuH^{s,q}(\Mbar,E) \to \cuH^{s-m,q}(\Mbar,E)
$$
for every $s \in \R$ and $1 < q < \infty$.
Here $\cuH^{0,q}(\Mbar,E) = {}^{\cu}L^q(\Mbar,E)$, the $L^q$-space of sections of $E$
with respect to a Hermitian metric on $E$ and the Riemannian density induced by a
cusp metric, i.e., a Riemannian metric ${}^{\cu}g$ in the interior of $\Mbar$ that in
a collar neighborhood of the boundary $Y$ takes the form
${}^{\cu}g = \frac{dx^2}{x^4} + g_Y$. For $m \in \N_0$ we have
$$
\cuH^{m,q}(\Mbar,E) = \{u \in {}^{\cu}L^q(\Mbar,E)\st Au \in {}^{\cu}L^q(\Mbar,E)
\textup{ for all } A \in \cuDiff^m(\Mbar;E)\},
$$
and for general $s \in \R$ the cusp Sobolev spaces $\cuH^{s,q}(\Mbar,E)$
are defined by duality and interpolation.

\subsection*{Cusp pseudodifferential operators}

By $\cuPsi^{\mu;\ell}(\Mbar,\Lambda)$ we denote the space of parameter-dependent families
of cusp pseudodifferential operators
$$
A(\lambda) : \dot{C}^{\infty}(\Mbar,E) \to \dot{C}^{\infty}(\Mbar,E)
$$
of order $\mu \in \R$, where $\Lambda$ is as usual a closed sector in $\C$.
This operator class can be described as follows: Let
$\omega,\tilde{\omega} \in C_c^{\infty}([0,\eps))$ be cut-off functions, i.e.,
$\omega,\tilde{\omega} \equiv 1$ near $x = 0$. We consider $\omega$ and $\tilde{\omega}$
functions on $\Mbar$ that are supported in the collar neighborhood
$U \cong [0,\eps)\times Y$ of the boundary.
\begin{itemize}
\item Whenever $\omega$ and $(1-\tilde{\omega})$ have disjoint supports, the operator
families $\omega A(\lambda) (1-\tilde{\omega})$ and $(1-\tilde{\omega})A(\lambda)\omega$
are integral operators with kernels
$$
k(z,z',\lambda) \in
\S\bigl(\Lambda,\dot{C}^{\infty}(\Mbar_z\times\Mbar_{z'},E\boxtimes E^*)\bigr).
$$
\item In the collar neighborhood $U$ of the boundary, the bundle $E|_U$ is
isomorphic to $\pi^*E|_Y$, where $\pi : [0,\eps)\times Y \to Y$ is the projection
on the second factor. Hence sections of $E$ on $U$ can be interpreted as
functions of $x \in [0,\eps)$ taking values in sections of $E|_Y$ on $Y$, e.g.,
$$
\dot{C}^{\infty}(U,E) \cong \dot{C}^{\infty}\bigl([0,\eps),C^{\infty}(Y,E|_Y)\bigr).
$$
The operator family $\omega A(\lambda) \tilde{\omega} : \dot{C}_c^{\infty}(U,E) \to
\dot{C}_c^{\infty}(U,E)$ now is of the form
$$
u \mapsto \frac{1}{2\pi}\int_{\R}\int_0^{\eps}e^{i(1/y-1/x)\xi} a(x,\xi,\lambda)u(y)\,
\frac{dy}{y^2}\,d\xi
$$
for $u \in \dot{C}_c^{\infty}\bigl([0,\eps),C^{\infty}(Y,E|_Y)\bigr)$, where
$$
a(x,\xi,\lambda) \in C^{\infty}\bigl([0,\eps)_x,L_{\cl}^{\mu;\ell}(Y,\R_{\xi}\times\Lambda;E|_Y)\bigr).
$$
\item The operator family $(1-\omega)A(\lambda)(1-\tilde{\omega})$ belongs to
the class $L_{\cl}^{\mu;\ell}(2\Mbar,\Lambda)$ of classical parameter-dependent
pseudodifferential operators of order $\mu$ on the double $2\Mbar$ of $\Mbar$
(acting in sections of an extension of the bundle $E$ to the double).
\end{itemize}
Every $A(\lambda) \in \cuPsi^{\mu;\ell}(\Mbar,\Lambda)$ extends by continuity to a
family of continuous operators
$$
A(\lambda) : \cuH^{s,q}(\Mbar,E) \to \cuH^{s-\mu,q}(\Mbar,E)
$$
for every $s \in \R$ and $1 < q < \infty$. The class of parameter-dependent families of cusp
pseudodifferential operators forms an algebra filtered by order.

\begin{theorem}\label{ParametrixThm}
Let $A \in \cuDiff^{m}(\Mbar;E)$, $m > 0$, be cusp-elliptic with parameter in $\Lambda$.
We have $A - \lambda \in \cuPsi^{m;m}(\Mbar,\Lambda)$, and there exists a
parameter-dependent parametrix $P(\lambda) \in \cuPsi^{-m;m}(\Mbar,\Lambda)$, i.e.,
$$
(A-\lambda)P(\lambda) - 1,\; P(\lambda)(A-\lambda) - 1 \in \cuPsi^{-\infty}(\Mbar,\Lambda).
$$
\end{theorem}

Theorem~\ref{ParametrixThm} is an instance of a standard result in the calculus of
pseudodifferential operators. The structural result that it entails about
parameter-dependent parametrices of $A - \lambda$ is the key to
the proof of Theorem~\ref{RbddnessResolvents} which is given in the next section.

\section{Proof of Theorem~\ref{RbddnessResolvents}}\label{ProofTheorem}

\noindent
Every $R(\lambda) \in \cuPsi^{-\infty}(\Mbar,\Lambda)$ gives rise to
an operator function
$$
R(\lambda) \in \S\bigl(\Lambda,\L(\cuH^{s,q}(\Mbar,E))\bigr)
$$
for every $s \in \R$ and $1 < q < \infty$.
Consequently, by Theorem~\ref{ParametrixThm}, the operator family
$$
A - \lambda : \cuH^{m,q}(\Mbar,E) \to {}^{\cu}L^q(\Mbar,E)
$$
is invertible for $\lambda \in \Lambda$ with $|\lambda| \geq R$ sufficiently large,
and, moreover,
$$
(A-\lambda)^{-1} - P(\lambda) : {}^{\cu}L^q(\Mbar,E) \to {}^{\cu}L^q(\Mbar,E)
$$
belongs to $\S\bigl(\Lambda_R,\L({}^{\cu}L^q(\Mbar,E))\bigr)$, where
$\Lambda_R = \{\lambda \in \Lambda \st |\lambda| \geq R\}$.
Here we are making use of the fact that the parameter-dependent parametrix
$P(\lambda)$ is tempered as a function of $\lambda \in \Lambda$ taking
values in the bounded operators on ${}^{\cu}L^q(\Mbar,E)$.

By Corollary~\ref{BilderRbounded} we get that
$$
\{\lambda(A-\lambda)^{-1} - \lambda P(\lambda) \st \lambda \in \Lambda_R\} \subset
\L\bigl({}^{\cu}L^q(\Mbar,E)\bigr)
$$
is $\RR$-bounded. To complete the proof it thus remains to show that
$$
\{\lambda P(\lambda) \st \lambda \in \Lambda_R\} \subset
\L\bigl({}^{\cu}L^q(\Mbar,E)\bigr)
$$
is $\RR$-bounded. To see this, let $\hat{\omega},\omega,\tilde{\omega} \in
C_c^{\infty}([0,\eps))$ be cut-off functions, i.e., $\hat{\omega},\omega,\tilde{\omega} 
\equiv 1$ near $x = 0$, and suppose that $\omega \equiv 1$ in a neighborhood of the
support of $\hat{\omega}$, and that $\tilde{\omega} \equiv 1$ in a neighborhood of the
support of $\omega$. We consider $\hat{\omega},\omega,\tilde{\omega}$ as functions on
$\Mbar$ that are supported in the collar neighborhood $U \cong [0,\eps)\times Y$ of
the boundary. Write
$$
P(\lambda) = \omega P(\lambda) \tilde{\omega} + (1-\omega)P(\lambda)(1-\hat{\omega})
+ R(\lambda).
$$
The operator family $R(\lambda) \in \cuPsi^{-\infty}(\Mbar,\Lambda)$, and so the set
$$
\{\lambda R(\lambda) \st \lambda \in \Lambda\} \subset \L\bigl({}^{\cu}L^q(\Mbar,E)\bigr)
$$
is $\RR$-bounded in view of Corollary~\ref{BilderRbounded}.
The family $(1-\omega)P(\lambda)(1-\hat{\omega})$ can be regarded as an element of
$L^{-m;m}_{\cl}(2\Mbar,\Lambda)$ (supported in the interior of the original copy
of $\Mbar$), and thus furnishes an element in
$S^{-m;m}_{\RR}\bigl(\Lambda;{}^{\cu}L^q(\Mbar,E),{}^{\cu}L^q(\Mbar,E)\bigr)$
in view of Corollary~\ref{PseudosClosedRbdd}. Consequently,
$$
\{\lambda(1-\omega)P(\lambda)(1-\hat{\omega}) \st \lambda \in \Lambda\} \subset
\L\bigl({}^{\cu}L^q(\Mbar,E)\bigr)
$$
is $\RR$-bounded. So it remains to show that
\begin{equation}\label{omegaPomegaRbdd}
\{\lambda \omega P(\lambda) \tilde{\omega} \st \lambda \in \Lambda\} \subset
\L\bigl({}^{\cu}L^q(\Mbar,E)\bigr)
\end{equation}
is $\RR$-bounded.

Consider $\omega P(\lambda) \tilde{\omega} : {}^{\cu}L^q(U,E) \to {}^{\cu}L^q(U,E)$,
where $U \cong [0,\eps)_x \times Y$.
Under the change of variables $r = -1/x$, this operator family can by definition of the
cusp calculus be regarded as the pull-back of an operator family
$$
\begin{gathered}
Q(\lambda) : L^q\bigl(\R_r,L^q(Y,E|_Y)\bigr) \to L^q\bigl(\R_r,L^q(Y,E|_Y)\bigr) \\
[Q(\lambda)u](r) = \frac{1}{2\pi}\iint e^{i(r-r')\varrho}a(r,\varrho,\lambda)u(r')
\,dr'\,d\varrho \quad \textup{for } u \in \S\bigl(\R_r,C^{\infty}(Y,E|_Y)\bigr),
\end{gathered}
$$
where $a(r,\varrho,\lambda) \in
S^0_{\cl}\bigl(\R_r,L^{-m;m}_{\cl}(Y,\R_{\varrho}\times\Lambda;E|_Y)\bigr)$.
In view of Corollary~\ref{PseudosClosedRbdd} we have
$$
L^{-m;m}_{\cl}(Y,\R\times\Lambda;E|_Y) \hookrightarrow
S^{-m;m}_{\RR}\bigl(\R\times\Lambda;L^q(Y,E|_Y),L^q(Y,E|_Y)\bigr),
$$
and so
$$
a(r,\varrho,\lambda) \in
S^0_{\cl}\bigl(\R_r,S^{-m;m}_{\RR}\bigl(\R_{\varrho}\times\Lambda;L^q(Y,E|_Y),L^q(Y,E|_Y)\bigr)\bigr).
$$
Thus, by Theorem~\ref{IterationRboundedness},
$$
Q(\lambda) = \op_r(a)(\lambda) \in
S^{-m;m}_{\RR}\bigl(\Lambda;L^q\bigl(\R,L^q(Y,E|_Y)\bigr),L^q\bigl(\R,L^q(Y,E|_Y)\bigr)\bigr),
$$
which shows that
$$
\{\lambda Q(\lambda) \st \lambda \in \Lambda\} \subset
\L\bigl(L^q\bigl(\R,L^q(Y,E|_Y)\bigr)\bigr)
$$
is $\RR$-bounded. Consequently, the set \eqref{omegaPomegaRbdd} is $\RR$-bounded, and
the proof of Theorem~\ref{RbddnessResolvents} is complete.

\end{document}